\documentclass[12pt,a4paper,twoside]{amsart}
\usepackage{amscd}
\usepackage{amssymb,bm}
\usepackage[scr=boondoxo]{mathalfa}
\usepackage{undertilde}

\theoremstyle{definition}

\newtheorem{defn}{Definition}[section]

\theoremstyle{remark}

\usepackage{amsfonts}

\renewcommand{\d}{\delta}
\newcommand{\D}{\Delta}

\DeclareMathOperator{\im}{im}

\usepackage{amsmath}
\usepackage{amssymb}
\usepackage{amsfonts}
\usepackage{enumerate}

\numberwithin{equation}{section}

\date{}

\def\char{\text{\rm char}}

\def\SL{\text{\rm SL}}
\def\GL{\text{\rm GL}}

\def\Ga{\Gamma}
\def\bbr{\mathbb{F}}
\def\bbf{\mathbb{R}}
\def\bbp{\mathbb{P}}

\def\bbn{\mathbb{N}}

\def\la{\lambda}

\def\suml{\sum\limits}
\theoremstyle{plain}
\newtheorem{thm}{Theorem}[section]
\newtheorem*{thm*}{Theorem}
\newtheorem{prop}[thm]{Proposition}
\newtheorem*{prop*}{Proposition}
\newtheorem*{prop**}{\ }
\def\beq{\begin{equation}}

  \def\ee{\end{equation}}

\theoremstyle{definition} 
\newtheorem{definition}[thm]{Definition}
\newtheorem*{definition*}{Definition}
\newtheorem*{thm1*}{Theorem A1}
\newtheorem*{thm2*}{Theorem A2}

\newtheorem*{conjecture*}{Conjecture}

\newtheorem*{claim*}{Claim}

\newtheorem*{remark*}{Remark}

\def\bbq{\mathbb{Q}}
\def\bbf{\mathbb{F}}
\def\bbr{\mathbb{R}}

\def\be{\begin{equation}}
\def\ee{\end{equation}}

\def\calh{\mathcal{H}}

\def\calm{\mathcal{M}}
\def\calb{\mathcal{B}}

\def\vare{\varepsilon}

\theoremstyle{remark}  
\begin{document}

\title[High Dimensional Expanders]{High Dimensional Expanders}
\author[A. Lubotzky]{Alexander Lubotzky}
\maketitle


\baselineskip 16pt

\begin{abstract}

Expander graphs have been, during the last five decades, the subject of a most fruitful interaction between pure mathematics and computer science, with influence and applications going both ways (cf. \cite{Lub94}, \cite{HLW06}, \cite{Lub12} and the references therein).  In the last decade, a theory of ``high dimensional expanders" has begun to emerge. The goal of the current paper is to describe some paths of this new area of study.

\end{abstract}

\section*{0. Introduction}

Expander graphs are graphs which are, at the same time,  sparse and highly connected.  These two seemingly contradicting properties  are what makes this theory non trivial and useful. The existence of such graphs is not a completely trivial issue, but by now there are many methods to show this: random methods, Kazhdan property $(T)$ from representation theory of semisimple Lie groups, Ramanujan  conjecture (as proved by Deligne and Drinfeld) from the theory of automorphic forms, the elementary Zig-Zag method  and ``interlacing polynomials".

The definition of expander graphs can be expressed in several different equivalent ways (combinatorial, spectral gap etc. - see \cite{Lub94}, \cite{Kam17a}).  When one comes to develop a high dimensional theory; i.e.\  a theory of finite simplicial complexes of dimension $d \ge 2$, which resembles that of expander graphs in dimension $d = 1$, the generalizations of the different properties are (usually) not equivalent.  One is led to notions like:  coboundary  expanders, cosystolic expanders, topological expanders, geometric expanders, spectral expanders etc. each of which has its importance and applications.

In \S 1, we recall very briefly several of the equivalent definitions of expander graphs (ignoring completely the wealth of their applications).  These will serve as pointers to the various high dimensional generalizations.

In \S 2, we will start with the spectral definition.  For this one needs ``discrete Hodge theory" as developed by Eckmann (\cite{Ec44}).  In this sense the classical work of Garland \cite{Gar73}, proving Serre's conjecture on the vanishing of the real cohomology groups of arithmetic lattices of $p$-adic Lie groups, can be considered as the earliest work on high dimensional expanders.  His ``local to global" method which treats the finite quotients of the Bruhat-Tits building has been rediscovered in recent years, with many applications, some of them will be described in \S 2.

In \S 3, we turn our attention to Gromov's topological and geometric expanders (a.k.a.~the  topological and geometric overlapping properties).   These quite intuitive directions were shown to be related to two much more abstract definitions of coboundary and cosystolic (high dimensional) expanders.  The last ones are defined using  the language of $\bbf_2$-cohomology.  Here also  a ``local to global" method enables to produce topological expanders from finite quotients of Bruhat-Tits buildings of $p$-adic Lie groups.

Section 4   will deal with random simplicial complexes, while in \S 5 we will briefly mention  several applications and connections with computer science.

\section{A few words about expander graphs}

Let $X  = (V,E) $ be a finite connected  graph with sets of vertices $V$ and edges $E$.  The Cheeger constant of $X$, denoted $h(X)$,  is:
$$ h(X) = \underset{A,B\subseteq V}{\inf} \; \frac{|E(A, B)|}{min(|A|, |B|)}$$
where the infimum  runs over all the possibilities of disjoint partitions $V = A \cup B$ and $E(A, B)$ is the set of edges connecting vertices in $A$ to vertices in $B$.

The graph $X$ is ${\pmb\vare}$-{\bf expander} if $h(X) \ge \vare$.

Let $L^2(V)$ be the space of real functions on $V$ with the inner product $\langle f, g\rangle = \suml_{v\in V} \deg (v) f(v) g(v)$  and $L^2_0 (V)$ the subspace of those which are orthogonal to the constant functions.  Similarly,  $L^2(E)$ is the space of functions on the edges with the standard inner product.

We fix an arbitrary orientation on the edges, and for $e \in E$ we denote its end points by $e^-$ and $e^+$.  Let $ d: L^2(V) \to L^2(E)$ be the map  $(df) (e) = f(e^+) - f(e^-) $ for  $f \in L^2(V)$ and $  \Delta = d^*d: L^2(V) \to L^2(V)$
when $d^*$ is the adjoint of $d$.  The operator $ \Delta$ is called the Laplacian of the graph.  One can show (cf. \cite[Chap.~4]{Lub94}), that it is independent of the chosen orientation.
One  can check that
$$ \Delta = I - M $$
when $M$ is the Markov operator on $L^2(V)$, i.e.,
$$ (Mf) (x) = \frac{1}{\deg(x)}\suml_{\{y|(x,y)\in E\}}f(y).$$

The smallest eigenvalue of $\D$ is 0 and it comes with multiplicity one if (and only if) $X$ is connected, which we will always assume.  The eigenfunctions with respect to 0 are the constant functions and as $\D$ is self adjoint, $L^2_0 (V)$ is invariant under $\D$  and the spectral gap
$$ \la_1(X) = \inf\left\{\frac{\langle\D f,f\rangle}{\langle f, f\rangle}\Big| f \in L^2_0 (V)\right\} $$
is the smallest eigenvalue of $\D$ acting on $L^2_0(V)$.

The following result is a discrete analogue of the classical Cheeger inequality (and its converse by Buser). This discrete version  was proved by Tanner, Alon and Milman (the reader is referred again to \cite[Chap.~4]{Lub94} for  a detailed history).

\begin{thm}\label{1.1} If $X$ is a finite connected $k$-regular  graph, then:
$$
\frac{h^2(X)}{2k^2} \le \la_1(X) \le \frac{2 h (X)}{k} $$
\end{thm}

We are usually interested in infinite families of $k$-regular graphs (``sparse").  Such a family forms a family of expanders (i.e., $h(X) \ge \vare$ for the same $\vare > 0$, for every $X$) if and only if $\la_1 (X) \ge \vare' > 0$ for the same $\vare'$ for every $X$.  I.e., Theorem \ref{1.1} says that expanders can be defined, equivalently, either by a combinatorial definition or using the spectral gap definition. Expressing this using the adjacency operator $A$ rather than the Laplacian $\D$: being expanders means that the second largest eigenvalue $\la(X)$ of $A$ is bounded away from $k$, which is the largest one.

Strictly speaking the notion of expanders requires spectral gap only in one side of the spectrum of $A$, but in many applications (e.g.\ if one wants to estimate the rate of convergence of the random walk on $X$ to the uniform distribution) one needs bounds on both sides.  Recall that $-k$ is also an eigenvalue of $A$ iff $X$ is bi-partite.  We can now define:  A $k$-regular connected graph is \emph{Ramanujan} if all  eigenvalues $\la$  of $A$ are either $\la = \pm k$ or $|\la| \le 2 \sqrt{k-1}$.  By the well-known Alon-Boppana theorem, the bound $2\sqrt{k-1}$ is the best one can hope for for an infinite family of $k$-regular graphs.  Let us recall that for the $k$-regular infinite tree $T = T_k$,  the classical result of Kesten asserts that the spectrum of the adjacency operator on $L^2(T_k)$ is exactly the interval $[-2\sqrt{k-1}, 2\sqrt{k-1}]$.
In a way Kesten's result lies beyond the Alon-Boppana theorem and there are many generalizations of this philosophy (cf. \cite{GrZu99}).

Ramanujan graphs were presented by Lubotzky-Phillips-Sarnak \cite{LPS88}, Margulis
\cite{Mar88}, Morgenstern and recently by Marcus-Spielman-Srivastava \cite{MSS15}.

There are several other ways to define expanders. Let us mention here one which has been observed only quite recently and has a  natural extension to high dimensional simplicial complexes.

Let $X$ be a finite connected $k$-regular graph, with adjacency matrix $A$, denote $k = q + 1$ and $$\mu(X) = \max\{ |\la|
\Big| \la \; \hbox{ e.v.\ of\ } A, \; \; \, \la \neq \pm k\}.$$
So $X$ is Ramanujan iff $\mu (X) \le 2\sqrt{q}$.  If $X$ is bipartite, write $V = V_0 \cup V_1$ where $V_0$ and $V_1$ are the two sides, and if not $V = V_0 = V_1$.  Let
$$ L^2_{00} (X) = \{ f \in L^2(V) \Big| \suml_{v\in V_i} f(v) = 0, \; \; \; \hbox{for\ } i = 0, 1\}.$$
So, $\mu(X)$ is the largest (in absolute value) eigenvalue of $A$ when acting on $L^2_{00} (X)$.
For $\la \in [ 2\sqrt{q}, q+1]$ write $\la = q^{1/p} + q^{(p-1)/p}$ for a unique $p\in [2, \infty]$, so $\lambda=2\sqrt{q}$ when $p=2$.

Now, let $\pi: T = T_k \to X$ be  a covering map.  For a fixed $t_0 \in T$, let $S_r = \{ t \in T\Big|distance (t,t_0) = r\}$ and for $f\in L^2_{00} (X)$ and $t\in T$, let $\tilde f(t) = \frac{1}{|S_r|} \suml_{s\in S_r} f (\pi(s))$ if $r = dist(t, t_0)$, i.e.\  $\tilde f$ is the averaging of the lift of $f$ around $t_0$.


\begin{thm}[Kamber \cite{Kam17a}]
$\mu (X) \le \la$ if and only if $\tilde f \in L^{p + \vare} (T)$ for all $f \in L^2_{00} (X), t_0 \in T$, and $\vare > 0$.
As a corollary $X$ is Ramanujan iff $$\tilde f \in L^{2 + \vare}(T),\, \forall t_0, \forall f, \forall \vare.$$
\end{thm}
\section{High dimensional expanders: spectral gap}

As described in \S 1, the notion of expander graphs can be expressed via a spectral gap property of the Laplacian.  This aspect has a natural high dimension version, but to present it one needs the language of real cohomology.  Let us start by recalling the basic notations.

\subsection{\label{sub:Simplicial-complexes-and-cohom}Simplicial complexes and
cohomology}

A finite simplicial complex $X$ is a finite collection of subsets, closed under inclusion,
of a finite set $X^{\left(0\right)}$, called the set of vertices of $X$.
 The sets in $X$ are called
\emph{simplices} or \emph{faces} and we denote by $X^{\left(i\right)}$
the set of simplices of $X$ of dimension $i$ ($i$-cells), which are the sets
in  $X$ of size $i+1$. So $X^{\left(-1\right)}$ is comprised of
the empty set, $X^{\left(0\right)}$ - of the vertices, $X^{\left(1\right)}$
- the edges, $X^{\left(2\right)}$ - the triangles, etc.
 Let $d = \dim X = \max \{ i |X^{(i)} \neq \emptyset\}$ and assume $X$ is a pure simplicial complex of dimension $d$, i.e., for every $F \in X$, there exists $G \in X^{(d)}$ with $F\subseteq G$.  Throughout
this discussion we will assume that $X^{\left(0\right)}=\left\{ v_{1},\ldots,v_{n}\right\} $
is the set of vertices and we fix an order $v_{1}<v_{2}<\ldots<v_{n}$
among the vertices. Now, if $F\in X^{\left(i\right)}$ we write $F=\left\{ v_{j_{0}},\ldots,v_{j_{i}}\right\} $
with $v_{j_{0}}<v_{j_{1}}<\ldots<v_{j_{i}}$. If $G\in X^{\left(i-1\right)}$,
we denote the \emph{oriented incidence number} $\left[F:G\right]$
by $\left(-1\right)^{\ell}$ if $F\backslash G=\left\{ v_{j_{\ell}}\right\} $
and $0$ if $G\nsubseteq F$. In particular, for every vertex $v\in X^{\left(0\right)}$
and for the unique face $\emptyset \in X^{\left(-1\right)}$, $\left[v:\emptyset\right]=1$.

If $\mathbb{F}$ is a field then $C^{i}\left(X,\mathbb{F}\right)$
is the $\mathbb{F}$-vector space of the functions from $X^{\left(i\right)}$
to $\mathbb{F}$. This is a vector space of dimension $\left|X^{\left(i\right)}\right|$
over $\mathbb{F}$ where the characteristic functions $\left\{ e_{F}\,\middle|\, F\in X^{\left(i\right)}\right\} $
serve as a basis.

The coboundary map $\delta_{i}:C^{i}\left(X,\mathbb{F}\right)\rightarrow C^{i+1}\left(X,\mathbb{F}\right)$
is given by:
\[
\left(\delta_{i}f\right)\left(F\right)=\sum_{G\in X^{\left(i\right)}}\left[F:G\right]f\left(G\right).
\]
So, if $f = e_{G}$ for some $G\in X^{\left(i\right)}$, $\delta_{i}e_{G}$
is a sum of all the simplices of dimension $i+1$ containing $G$
with signs $\pm1$ according to the relative orientations.

It is well known and easy to prove that $\delta_{i}\circ\delta_{i-1}=0$.
Thus $B^{i}\left(X,\mathbb{F}\right)=\im\delta_{i-1}$ - ``the space
of $i$-coboundaries'' is contained in $Z^{i}\left(X,\mathbb{F}\right)=\ker\delta_{i}$
- the $i$-cocycles and the quotient $H^{i}\left(X,\mathbb{F}\right)=
{Z^{i}\left(X,\mathbb{F}\right)}/{B^{i}\left(X,\mathbb{F}\right)}$
is the $i$-th cohomology group of $X$ over $\mathbb{F}$.

In a dual way one can look at $C_{i}\left(X,\mathbb{F}\right)$ -
the $\mathbb{F}$-vector space spanned by the simplices of dimension
$i$. Let $\partial_{i}:C_{i}\left(X,\mathbb{F}\right)\rightarrow C_{i-1}\left(X,\mathbb{F}\right)$
be the boundary map defined on the basis element $F$ by: $\partial F=\sum_{G\in X^{\left(i-1\right)}}\left[F:G\right]\cdot G$,
i.e.\ if $F=\left\{ v_{j_{0}},\ldots,v_{j_{i}}\right\} $ then $\partial_{i}F=\sum_{t=0}^{i}\left(-1\right)^{t}$\newline $\left\{ v_{j_{0}},\ldots,\widehat{v_{j_{t}}},\ldots,v_{j_{i}}\right\} $.
Again $\partial_{i}\circ\partial_{i+1}=0$ and so the boundaries $B_{i}\left(X,\mathbb{F}\right)=\im\partial_{i+1}$
are inside the cycles $Z_{i}\left(X,\mathbb{F}\right)=\ker\partial_{i}$
and $H_{i}\left(X,\mathbb{F}\right)=
{Z_{i}\left(X,\mathbb{F}\right)}/{B_{i}\left(X,\mathbb{F}\right)}$
gives the $i$-th homology group of $X$ over $\mathbb{F}$. As $\mathbb{F}$
is a field, it is not difficult in this case to show that $H_{i}\left(X,\mathbb{F}\right)\simeq H^{i}\left(X,\mathbb{F}\right)$.

In the next section, we will need the case $F = \bbf_2$ - the field of two elements, but for the rest of Section 2 we work with $F = \bbr$.
In this case $C^i(X,\bbr)$ has the natural structure of a Hilbert space, where for $f, g \in C^i(X,\bbr), \; \; \langle f, g\rangle = \suml_{F \in X^{(i)}} \deg (F) f (F) g (F)$, when $\deg(F) = \# \{ G \in X^{(d)} \Big | G \supseteq F\}$.  Now, $C_i(X,\bbr)$ is the dual of $C^i (X,\bbr)$ in a natural way and we can identify them and treat the operators $\D^{up}_i = \delta^*_i\delta_i, \D^{down}_i = \delta_{i-1} \delta^*_{i-1}$ and $\D_i = \D^{up}_i + \D^{down}_i$ as operators from $C^i$ to $C^i$, all are self-adjoint with non-negative eigenvalues. One may check that $$\left(\delta_{i}^{*}f\right)\left(G\right)=\frac{1}{\deg\left(G\right)}\sum_{F\in X^{\left(i+1\right)}}\left[F:G\right]\deg\left(F\right)f\left(F\right)$$ for $f\in C^{i+1}\left(X,\mathbb{R}\right)$ and $G\in X^{(i)}$, so in the regular case $\delta_i^*$ is equal to $\partial_{i+1}$ up to a constant multiple. Define $\mathcal{Z}_i=\ker\delta^*_{i-1}$ and $\mathcal{B}_i=\im\delta_i^*$ (so in the regular case $\mathcal{Z}_i=Z_i, \mathcal{B}_i=B_i$). The following proposition, going back to Eckmann \cite{Ec44}, is elementary:
\begin{prop}[Hodge decomposition] $C^i = B^i \oplus \mathcal{H}^i \oplus \mathcal{B}_i $ when $\calh^i = Ker (\D_i)$ is called the space of Harmonic cycles.  In fact $\calh^i \simeq H^i (X, \bbr)$.  Note that $\D^{up}_i$ vanishes on $\mathcal{Z}^i = B^i \oplus \mathcal{H}^i$.
\end{prop}
\begin{definition}
The {\em $i$-dimensional spectral gap} of $X$ is $\la^{(i)} (X) = \min\{ \la\Big| \la\; \hbox{ e.v. of\ } \D^{up}_i\Big|_{(B^i)^\perp}\}$.  One may check that $(B^i)^\perp = \mathcal{Z}_i$, and as $\D^{up}_i = \delta^*_i\circ \delta_i$, we have $$\la^{(i)} (X) =
\underset{f \in (B^i)^\perp}{\inf}\left\{
\frac{|\langle\D^{up}_i f, f\rangle|}{\langle f, f\rangle} \right\} =
\left(\underset{f \in \mathcal{Z}_i}{\inf} \left\{\frac{\| \delta f\|}{\| f \|}\right\}\right)^2.$$
Also, $\D^{up}_i$ vanishes also on $\calh^i$, so  $\la^{(i)} (X) > 0$
  implies $H^i (X, \bbr) = \{ 0 \}$, and the converse is also true.
\end{definition}

For a  $k$-regular graph $(B^0)^\perp = Z_0 = L^2_0 (X) $ and so $\la_1 (X)$ that was defined in \S 1 for a graph $X$, is
  $\la^{(0)} (X)$ in the notations here.  We define:
\begin{definition}
A pure $d$-dimensional simplicial complex will be called $\vare$-\emph{spectral expander} if for every $i = 0,\dots, d-1, \; \; \la^{(i)}\ge \vare$.
\end{definition}

Recall that the Alon-Boppana theorem asymptotically bounds the spectral gap of $k$-regular graphs by that of their universal cover, the $k$-regular tree. In higher dimension the situation is more involved:
\begin{thm}[\cite{PR17}] For an infinite complex $X$, let $\lambda^{(i)}(X)$ be the
bottom of the spectrum of $\Delta^{up}_i (X)$ on $(B^i)^\perp$.
       Let $\{X_n\}$ be a family of quotients of $X$, such that the injectivity radius of $X_n$ approaches infinity. If zero is not an isolated point in the spectrum of $\D_i^{up}(X)$, on $(G^i)^\perp$, then $$\liminf_{n\rightarrow\infty}\{\lambda^{(i)}(X_n)\}\le \lambda^{(i)}(X).$$
\end{thm}
Note that zero cannot be an isolated point in the spectrum of the Laplacian of an infinite graph, since the constant function is not in $L^2$. However, for complexes of higher dimensional this can happen, and in this case the Alon-Boppana principle can be violated (see \cite[Thm.\ 3.10]{PR17} for an example).

\subsection{Garland method}

The seminal paper of Howard Garland (\cite{Gar73}, see also \cite{Bor73}), can be considered as the first paper on high dimensional expanders.  It gave  examples of spectral expanders, by a method which bounds the eigenvalues of the simplicial complex by the eigenvalues of its links.  Garland's method has been revisited in recent years with various simplifications and extensions.  Let us give here one of them, but we need more definitions:
If $F$ is a face of $X$  of dimension $i$, the  \emph{link of $F$ in $X$} denoted $\ell k_X(F)$, is
\[ \ell k_X(F): = \{ G \in X \Big| F \cup G \in X, F\cap G = \emptyset\}.\]
One can easily check that if $X$ is a pure simplicial complex  of dimension $d$,
$\dim(\ell k_X (F)) = d - i - 1$.

Garland's method can be conveniently summarized by the following theorem.  Note that if $dim(X) = d$ and $dim (F) = d -2$, then $\ell k_X(F)$ is a graph.
\begin{thm}[Garland, cf. \cite{GuWa16} ] If $dim (X) = d$ and for every face $F$ of dimension   $d-2$, $\la^{(0)}(\ell k_X(F)) \ge \vare$, then
\[ \la^{(d-1)} (X) \ge 1+d\, \vare - d.\]
\end{thm}

So, Garland's method enables to give a fairly good bound on $\la^{(d-1)} (X)$ if all links of $d-2$ faces are very good expanders.  One can use the result to bound also $\la^{(j)} (X)$ for $j \le d-1$, by replacing $X$ with its $j+1$ skeleton, i.e., the collection of all the faces of $X$ of dimension at most $j+1$.
In fact, even more: if the links of the $(d-2)$-faces are excellent expander graphs and the 1-skeleton is connected, then the complex is spectral expander (cf. \cite{Op17b}).
In the next subsection, we will explain Garland's motivation and results.  But in recent years his method have been picked up in various different directions.  Most of them have to do with vanishing of some cohomology groups.

One of the nicest applications of Garland's method is the work of Zuk \cite{Zu03}, Pansu\cite{Pa89} and  Ballman-Swiatkowski \cite{BaSw97}.  The starting point of these works is the well-known result that a discrete group $\Ga$ has Kazhdan property $(T)$ iff $H^1(\Ga, V) = \{ 0 \}$ for every unitary representation of $\Ga$ on any Hilbert space.  These authors used Garland's work to deduce such a vanishing result for $H^1$ if $\Ga$ acts cocompactly on an infinite contractible  simplicial complex of dimension 2 all of whose vertex links are very good expanders.  The most amusing is Zuk's method which enables (sometimes) to deduce property $(T)$ from a presentation of $\Ga$ by generators and relations.   For example it shows property $(T)$ for some random groups (see also \cite{KK13}).  This is very different than the way Kazhdan  produced the first groups with property $(T)$ and it shows that property $(T)$ is not such a rare property.

A work of a similar flavor but in a different direction is the work of De Chiffre,  Glebsky, Lubotzky, and Thom (\cite{DGLT}).  Recall first (vaguely) the basic definition of ``group stability": Consider the degree $n$ unitary group $U(n)$ with an invariant metric $d_n$.  We say  that a group $\Ga$ presented by  a finite set of generators   $S$ with finitely many relations $R$, is $\left(U(n), d_n\right)$-stable if every  \emph{almost representation} $\rho$ of $\Ga$ into $U(n)$  is close to a representation $\tilde \rho$.  By ``almost"  we mean that $\rho(r)$ is very close to the identity for every $r \in R$ and ``close" means that $\rho(s)$ and $\tilde \rho(s)$ are close w.r.t. $d_n$, for every $s \in S$.   One can study these questions w.r.t. different distance functions, e.g., the one induced by the Hilbert-Schmidt norm, the operator norm or the $L^2$-norm, a.k.a.~the Frobenius norm.

Let us stick to the $L^2$-norm. In  \cite{DGLT}  it is shown that if $H^2(\Ga, V) = \{ 0 \}$ for every unitary representation of $\Ga$, then $\Ga$ is $\left(U(n), d_{L^2}\right)$-stable.  Then the Garland method is used (along the line of the results mentioned above for $H^1$) to produce many examples of $L^2$-stable groups by considering actions on 3-dimensional infinite simplicial complexes, whose edge-links are excellent expanders. This implies that many high rank  cocompact lattices in simple $p$-adic Lie groups are $(U(n), d_{L^2})$-stable. The most striking application is proving that there exists  a group which is not $L^2$-approximated (the reader is referred to \cite{DGLT}  for the definitions and exact results and to \cite{T18} for background and applications).

In \cite{GuWa16}, Gundert and Wagner used the Garland method to estimate the eigenvalues of random simplicial complexes - see also \S 4.
For some stronger versions of Garland's method - see \cite{Op17a}, \cite{Op17b} and the references therein.

\subsection{Bruhat-Tits buildings and their finite quotients}

Let $K$ be a non-Archimedean local field, i.e., $K$ is a finite extension of $\bbq_p$, the field of $p$-adic numbers, or $K$ is $\bbf_q((t))$-the field of Laurent power series over a finite field $\bbf_q$. Let $\mathcal{O}$ be the ring of integers of $K, \, \calm$ the (unique) maximal ideal of $\mathcal{O}$, and $\bbf_q = \mathcal{O}/\calm$ the finite quotient  where $q = p^\ell$ for some prime $p$ and $\ell\in\bbn$. Let $\utilde G$
be a $K$-simple simply connected group of $K$-rank $r$, e.g., $\utilde G = \SL_n$ in which case $r = n - 1$, and let $G = \utilde G(K)$.

Bruhat and Tits developed a theory which associates with $G$ an infinite (if $r \ge 1$) contractible simplicial complex $\mathcal{B} = \mathcal{B}(G)$ of dimension $r$.  Here is a quick description of it:  $G$ has $r+1$ conjugacy classes of maximal compact subgroups (cf. \cite[Theorem~3.13, p.~150]{PlRa91}) and a unique class of maximal open pro-$p$ subgroups, called Iwahori subgroups.  The vertices of $\mathcal{B}$ are the maximal compact subgroups (so they come with $r+1$ ``colors" according to their conjugacy class) and a set of $i + 1$ such vertices form a cell if their intersection contains an Iwahori subgroup.  This is an $r$-dimensional simplicial complex whose maximal faces can be identified with $G/I$ when $I$ is a fixed Iwahori subgroup (for more see \cite{BrTi72}, \cite{PlRa91} and \cite{Lub14} for a quick explicit description of $\mathcal{B}\left(\SL_n(\bbq_p)\right)$.  The case of $\mathcal{B} \left(\SL_2(K)\right)$, which is a $(q+1)$-regular tree, is studied in detail in \cite{Se80}).

Let $\Ga$ be a cocompact lattice in $G$, i.e., a discrete subgroup with $\Ga\setminus G$ compact.  Assume, for simplicity that $\Ga$ is torsion free, a condition which can always be achieved by passing to a finite index subgroup.  Such $\Ga$ is always an arithmetic lattice if $r \ge 2$ by Margulis arithmeticity Theorem (\cite{Mar91}) and, at least if $\char(K)=0$, there are always such lattices by Borel and Harder (\cite{BoHa78}).  When we fix $K$ and $G$ and run over all such lattices in $G$, for example, over the infinitely many congruence subgroups of $\Ga$, we obtained a family of \emph{bounded degree} simplicial complexes, i.e.\  every vertex is included in a bounded number of faces.  These simplicial complexes, give the major examples of ``high dimensional expanders" discussed in this paper.

Garland's method described in the previous subsection was developed by him in order to prove a conjecture of Serre asserting that if $r \ge 2$, $H^i(\Ga, \bbr) = \{ 0 \}$ for every $\Ga$ as above and every $1 \le i \le r -1$.  Indeed, the vertex links of the building $\mathcal B$ are the associated spherical building over the finite field $\bbf_q$ (for example, for $\utilde G = \SL_n$, this is the flag complex of the proper subspaces of $\bbf^n_q$). For such buildings, for every cell $F\in X^{(i)}, 0\le i\le d-2$, one has $\la^{(0)} \left(\ell k(F)\right) \to 1$ when $q \to \infty$ (e.g.\ for $\utilde G = \SL_3$, we get the $(q+1)$-regular ``points to lines graph" of the projective plane $\bbp(\bbf^3_q)$, for which one can check that $\la_1 = 1 - \frac{1}{\sqrt{q}}$. See \cite{Gar73, BaSw97, EvKa16}).  One therefore can deduce from Theorem 2.5 that if $q\ge q(\utilde G)$, then $\la^{(i)} (X) > \vare'$ for every $i = 1, \dots, r-1$ and every finite quotient $X$ of $\mathcal B = \mathcal{B} (G)$. In particular, all these quotients are spectral expanders as defined in Definition 2.3.

This also implies Serre's conjecture if $q$ is sufficiently large (see Definition 2.2).  Serre's conjecture has been proved in full since then (cf.~\cite{Ca74}
 and \cite[Chap.~XI]{BoWa80}) by representation theoretic methods, but Garland's method has its own life in various other contexts.

In Section 1, Theorem 1.2,  we saw that expander graphs can also be defined as ``$L_p$-expanders" for a suitable $2 \le p \in \bbr$.  This definition can be extended to high dimensional simplicial complexes and is especially suitable in the context of this subsection.

Let $\mathcal{B}$ be one of the Bruhat-Tits buildings described above and $\pi: \mathcal{B} \to X$ the covering map. Let $f \in L^2_0 (X^{(r)})$, i.e.\  a function orthogonal to the constants on the $r$-cells of $X$ (one can consider also $i$-cells for $0 \le i \le r$, but we stick to these for simplicity of the exposition, the reader is referred to \cite{Kam17b} for a more general setting).  Now, using the notion of $W$-distance on $\mathcal{B}^{(r)}$, when $W$ is the affine Weyl group of $G$, one can define for a fixed $t_0 \in \mathcal{B}^{(r)}$, a function $\tilde f $ on $\mathcal{B}^{(r)}$ - the $r$-faces of $\mathcal{B}$, as follows:  For $t\in\mathcal{B}^{(r)}$, let $\tilde f (t) = \frac{1}{|S_t|} \sum\limits_{s \in S_t} f\left(\pi(s)\right)$  when $S_\ell = \{ s \in \mathcal{B}^{(r)} \Big| W\hbox{-distance\ }  (s, t_0 ) = W\hbox{-distance\ } (t, t_0)\}.$

\begin{definition}\label{Oh} We say that $X$ is $L_p$-expander if for every $t_0$ and $f$ as above $\tilde f \in L^{p+\vare} (\mathcal {B}^{(r)})$ for every $\vare > 0$.

\end{definition}
Applying Oh's result \cite{Oh02} which gives the exact rate of decay of the matrix coefficients of the unitary representations of $G$, the so-called ``quantitative property $(T)$", Kamber deduced that $X$ as above are always $L_p$-expanders when $p = p(G)$ according to the following table.

\bigskip

\begin{center}
\begin{tabular} { |c|c|c|c|c|c|c|c|c|c|c| }
 \hline
 $W$ &$\tilde A_n$ &$\tilde B_n$ &$\tilde C_n$ &$\tilde D_n, \, n \; even$ &$\tilde D_n,\,  n \; odd$ &$\tilde E_6$ &$\tilde E_7$ &$\tilde E_8$ &$\tilde F_4$ &$\tilde G_2$\\
 \hline
  $p$ &$2n$ &$2n$ &$2n$ &$2(n-1)$  &$2n$ &$16$ &$18$ &$29$ &$11$ &$6$ \\
\hline
\end{tabular}
\end{center}


\bigskip
Let us stress that this is not just an abstract result.  From this fact, we can deduce non-trivial inequalities on the eigenvalues of various ``Hecke operators" acting on the faces of $X$.  The reader is referred to \cite{Kam17b} for more in this direction.

\subsection{Ramanujan complexes}

Ramanujan graphs stand out among expander graphs as the optimal expanders from a spectral point of view (cf. \cite{Val97}).  These are the finite connected $k$-regular  graphs $X$ for which every eigenvalue $\la$ of the adjacency matrix $A = A_X$ satisfies either $|\la| = k $ or $|\la | \le 2\sqrt{k-1}$.  The first constructions of such graphs were presented as an application of the works of Deligne (in characteristic zero) and Drinfeld (in positive characteristic) proving the so called Ramanujan conjecture for $\GL_2 $ (see \cite{Lub94} for a detailed survey).  Recently, a new (non-constructive) method has been presented in \cite{MSS15}.

It is therefore not surprising that following the work of Laurent Lafforgue \cite{Laf02} (for which he got the Fields Medal) extending Drinfeld's work from $\GL_2 $ to $\GL_d$, general $d$, several mathematicians have started to develop a high dimensional theory of Ramanujan simplicial complexes, cf. (\cite{CSZ03}, \cite{Li04}, \cite{LSV05a}, \cite{LSV05b}, \cite{Sar07}).  One may argue what is ``the right" definition of Ramanujan complexes (see the above references and \cite{KLW10}, \cite{Kan16}, \cite{Fi16}, \cite{Kam17b}, \cite{LLP17}).  This topic deserves a survey of its own.  Here we just briefly point out  some directions of research which came out in the work of several mathematicians.

In the context of $X=\Ga\setminus\calb$ where $\calb$ a Bruhat-Tits building associated with $G =\utilde G (K)$ as in \S 2.3,  and $\Ga$ a cocompact lattice acting on it, the most sensible definition seems to be the following:

\begin{definition} In the notation above, $\Ga\setminus\calb$ is called a Ramanujan complex if every infinite dimensional irreducible $I$-spherical $G$-subrepresentation  of $L^2(\Ga\setminus G)$ is tempered.
\end{definition}

Recall that $I$ is the Iwahori subgroup defined above, a representation is $I$-spherical if it contains a non-zero $I$-fixed vector  and it is tempered if it is weakly-contained in $L^2(G)$.

This definition can be expressed also in other ways; it is $L_2$-expander in the notations of \cite{Kam17b} and Definition \ref{Oh} above.  It can also be expressed in a combinatorial-spectral way.  For the group $\SL_2$, in which case $\calb$ is a $(q+1)$-regular tree and $X = \Ga\setminus\calb$ is a $(q+1)$-regular graph, this definition is equivalent to the graph being Ramanujan graph.  Ramanujan complexes are also optimal among high-dimensional expanders (see \cite{Li04}, \cite{LSV05a} and \cite{PR17}).  For most applications so far (such as the geometric and topological expanders to be presented in \S 3) one does not need the full power of the Ramanujan property and quantitative Property $(T)$ (\'a la Oh \cite{Oh02}, see \S2.3) suffices. On the other hand the study of the cut-off phenomenon of Ramanujan complexes in \cite{LLP17} did use the full power of the Ramanujan property.  The same can be said about the application of Ramanujan graphs and Ramanujan complexes to the study of ``golden gates" for quantum computation (see \cite{PS17}  and \cite{PS18}), where the Ramanujan bounds give a distribution of elements in $SU(2)$ with ``optimal entropy".

The Ramanujan graphs of \cite{LPS88} (a.k.a.~the LPS-graphs) have also been used to solve other combinatorial problems.  For example they give the best (from a quantitative point of view) known examples of ``high girth, high chromatic number" graphs.  After finding the appropriate high dimensional notions of ``girth" and ``chromatic number", these results can indeed be generalized to the Ramanujan complexes constructed in \cite{LSV05b}, (see \cite{LM07}, \cite{GP14}, \cite{EGL15}).

Ramanujan graphs can be characterized as those graphs whose associated zeta functions satisfy ``the Riemann Hypothesis (RH)" - see \cite{Lub94},  for an exact formulation and references.  An interesting direction of research is to try to associate to high dimensional complexes suitable ``zeta functions" with the hope that also in this context the Ramanujaness of the complex can be expressed via the RH.  For this direction or research - see \cite{Sto06}, \cite{KaLi14}, \cite{DK14}, \cite{KLW10}, \cite{Kan16}, \cite{Kam17b} and \cite{LLP17}.

\section{Geometric and Topological expanders}

In this chapter we will describe a phenomenon which is truly high dimensional; the geometric and topological overlapping properties which lead to geometric and topological expanders.  The latter call for coboundary and cosystolic expanders.

\subsection{Geometric and  Topological overlapping}
Our story begins with a result of Boros  and F\"uredi \cite{BF84}, at the time two undergraduates in Hungary, who proved the following result, as a response to a question of Erd\"os:  If $P$ is a set of $n$ points in $\bbr^2$, then there exists a point $z \in \bbr^2$ which is covered by $\bigg(\frac 2 9 - o(1)\bigg) {n\choose 3}$ of the ${n\choose 3}$ affine triangles determined by these points.  Shortly afterward B\'arany \cite{Bar82} proved the $d$-dimensional version:  For every $d \in \bbn, \exists\;  0 < C_d \in \bbr$, such that if $P \subseteq \bbr^d$ with $|P|= n$, then there  exists $z \in\bbr^d$ which is covered by at least $C_d{n \choose {d+1}} $ of the ${n\choose {d+1}}$ affine simplices determined by these points.

While $2/9$ is optimal for $d = 2$, it is not known what are the optimal $C_d$'s, neither what is their rate of convergence to 0, when $d$ goes to infinity.

B\'arany's result can be rephrased as:  Let $\D^{(d)}_n$ be the complete $d$-dimensional simplicial complex on $n$ vertices (i.e.\  the collections of all subsets of $[n]$ of size at most $d+1$) and $f: \D^{(d)}_n \to \bbr^d$ an affine map. Then there exists $z \in\bbr^d$ which is covered by at least $C_d {n\choose d+1}$ of the images of the $d$-dimensional faces.

In \cite{Gro10}, Gromov proved the following amazing result:  B\'arany's theorem above is true for every continuous map $f: \D^{(d)}_n \to \bbr^d$.  In fact, he proved it with constants $C_d \in \bbr$ which were better than what was known before for affine maps.  The reader is encouraged to draw the 2-dimensional case to realize how surprising and even counter-intuitive this theorem is!  Gromov also changed the point of view on these types of results; rather than thinking of them as  properties  of $\bbr^d$, think of them as properties of the simplicial complex $X$.  Let us now define:

\begin{definition} A $d$-dimensional pure simplicial complex $X$ is said to be $\vare$-\emph{geometric}  (resp. $\vare$-\emph{topological}) \emph{expander} if for every affine (resp. continuous) map $f : X \to \bbr^d$, there exists $z \in\bbr^d$ such that $\vare$-proportion of the images of the $d$-cells in $X^{(d)}$, covers the point $z$.
\end{definition}
So B\'arany (resp., Gromov) Theorem is the claim that $\D^{(d)}_n$, the complete simplicial complex of dimension $d$ on $n$ vertices, is $C_d$-geometric (resp., $C_d$-topological) expander.

Let us look for a moment at the case of dimension one to see why we call this property ``expander":  If $X = (V, E)$ is an expander graph and $f:X\to\bbr$ any continuous map, choose a point $z \in \bbr$ such that the two disjoint sets
\[ A = \{ v \in V\Big| f(v) < z\} \; \hbox{ and\ } B = \{ v \in V\Big| f(v) > z\} \]
are of size approximately $\frac{|V|}{2}$.  By the expansion property, there are many edges in $E$ which connect $A$ and $B$. The image of each such an edge under $f$ must pass through $z$ by the mean value theorem.  Hence $X$ is also a topological expander.

      We should mention that a topological expander graph $X$ does not have to be an expander graph.  Moreover,  it does not even have to be  connected.  For example, assume $X$ is a union of a large expander graph and another small (say of size $o(|X|)$) connected component.  Then $X$ is a topological expander even though it is not an expander graph.

Anyway, Gromov and B\'arany Theorems refer to the complete simplicial complexes:  note how difficult is the case $d \ge 2 $ and how trivial it is to prove that the complete graph is an expander.  Gromov also proved that some other interesting simplicial complexes are $d$-dimensional topological expanders, e.g., the flag complexes of $d+2$ dimensional vector spaces  over finite fields or more generally spherical buildings of simple algebraic groups over finite fields (cf. \cite{LMM16}).  All these examples are \emph{not} of bounded degree.  Recall (see also Definition 4.1 below)  that we say that a family of $d$-dimensional simplicial complexes are of \emph{bounded degree} (resp. \emph{bounded upper degree}) if  for every vertex $v$ (resp., every face $F$ of dimension $d-1$) the number of faces containing it is bounded.  The non trivial aspect of expander graphs in dimension one is the construction of such graphs of bounded degree.

 Gromov \cite{Gro10}  put forward the basic questions:  Let $d \ge 2$, are there bounded degree $d$-dimensional geometric/topological expanders?

 The existence of geometric expanders of bounded degree was shown by Fox, Gromov, Lafforgue, Naor and Pach \cite{FGLNP12} in several ways - most notably are two:  the random method which we will come back to in \S 4 and the second is by showing that for a fixed $d$, if $q$ is a large enough (depending on $d$) and fixed,  the Ramanujan complexes described in \S 2.4  are geometric expanders of bounded degree.  A more general version was given by Evra \cite{Ev17} .
 \begin{thm}\label{th3.2} Given $2\le d \in \bbn$,  there exists $q_0 = q_0(d)$ and $\vare = \vare(d)$ such that for every  $q > q_0$, if $K$ is a non-Archimedean local field of residue degree $q$ and $\utilde G$ a simple $K$-algebraic group of $K$-rank $d$, then the finite quotients of $\mathcal{B} = \mathcal{B} (\utilde G(K))$ - the Bruhat-Tits building associated with $G = \utilde G(K)$ - are all $\vare$-geometric expanders.
 \end{thm}

Theorem 3.2 is deduced in \cite{Ev17} in a similar way as the proof in \cite{FGLNP12} using a ``mixing lemma" and a classical convexity result of Pach \cite{Pac98}.  The mixing lemma is deduced there from Oh's ``quantitative property $(T)$ (\cite{Oh02}).   The language of $L_p$-expanders described in \S 2 gives a systematic way to express this. (Compare also to \cite{PRT16}). The fact that we have an $\vare = \vare(d)$ which is independent of $q$, provided $q > q_0$, (which is more than one needs in order to answer Gromov's geometric question) is due to the fact that for a fixed  $d \in \bbn$, one has the same $p$ in the table in \S 2.3.  which works for all groups of rank $d$.

The question of bounded degree topological expanders is much more difficult and will be discussed in the next subsections.

\subsection{Coboundary expanders}

As of now there is only one known method (with several small variants) to prove that a simplicial complex $X$ is a  topological expander.  This is via ``coboundary expander" which requires the language of cohomology as introduced in \S 2.1, but this time with $\bbf_2$-coefficients.

Let $X$ be a finite $d$-dimensional pure simplicial complex, define on it a \emph {weight function} $w$ as follows: for $F \in X^{(i)}$ let \[ w(F) = \frac{1}{{d+1\choose i + 1} |X^{(d)}|} \left| \{ G \in X^{(d)} | G \supseteq F\} \right|.\]
One could work with a number of different weight functions, but this one is quite pleasant, for example, it is a probability measure on $X^{(i)}$; one easily checks that $\suml_{F \in X^{(i)}} w(F) = 1$.  Now for $f \in C^i(X, \bbf_2)$,  denote $\| f \| =  \mathop{\sum w(F)}\limits_{\{F \in X^{(i)}| f(F) \neq 0\}}.$
We can now define the important notion of ``coboundary expanders" - a notion which was independently defined by Linial-Meshulam \cite{LM06} and Gromov \cite{Gro10} (in both cases without calling it coboundary expanders).

\begin{definition} Let $X$ be as above:\\

\begin{enumerate}[(a)]

 \item  For $0 \le i \le d-1$, define the $i$th coboundary expansion  $\mathscr{h}_i(X)$ of $X$ as:
\[ \mathscr{h}_i(X) = \min\limits_{f \in C^i\setminus B^i} \frac{\| \delta_i f\|}{\| [f]\|}\]
where $[f] = f + B^i$ is the coset of $f$ w.r.t. the $i$-coboundaries  and $\| [f]\| = \min\limits_{g \in [f]} \| g \|$.
(Note that $\| [f]\|$ is the ``normalized distance" of $f$ from $B^i$).
Let $\mathscr{h}(X) = \min\{\mathscr{h}_i(X) | i = 0, \dots, d-1\}$.

\item The complex $X$ is said to be $\vare$\emph{-coboundary expander} if $\mathscr{h} (X) \ge \vare$.
\end{enumerate}
\end{definition}

A few remarks are in order here:

\begin{enumerate}[(i)]

\item The reader can easily check that if $X$ is a $k$-regular graph, then
$ \mathscr{h}(X) =\frac 2k \cdot h(X)$ where $h(X)$ is the Cheeger constant of the graph as defined in \S 1.  So, indeed the above definition extends the notion of expander graphs.

    \item The definition of $\vare_i$, and especially the fact that the minimum runs over $f \in C^i\setminus B^i$ looks unnatural at first sight, but if we recall that $\| [f]\|$ is exactly the ``norm" of the element in $f + B^i$ which is closest to $B^i$, we see that this corresponds to going over $(B^i)^\perp$ when we consider real coefficients.  Moreover as pointed out in \S 2, over $\bbr$,
        \[\frac{\langle \D^{up}_i f, f\rangle}{\langle f, f\rangle} = \frac {\| \delta_if\|}{\| f\|}\] and so the definition of $\mathscr{h}_i$ here is
``the characteristic 2 analogue" of the spectral gap defined in Definition 2.2.  For the connection between the spectral gap and the coboundary expansion - see \cite{SKM14}, \cite{PRT16} and \cite{GS15}.

\item Also here it is easy to see that $\mathscr{h}_i (X) > 0 \; \hbox{ iff\ } H^i(X, \bbf_2) = \{ 0\}$.
\end{enumerate}

A basic result proved independently in \cite{LM06}, \cite{MW09} and \cite{Gro10} is:
\begin{thm} for the complete $d$-dimension complex $\D^{(d)}_n, \mathscr{h}_i (\D^{(d)}_n) \ge 1 - o_d(1)$ for every $i = 0, \dots, d-1$.
\end{thm}

Here is Gromov fundamental result on the connection between coboundary expanders and topological expanders:
\begin{thm} Coboundary expanders are topological expanders, namely,  for every $d \in\bbn$ and $0 < \vare \in\bbr$, there exists $\vare^1 = \vare^1(d, \vare) > 0$ such that if $X$ is a $d$-dimensional complex which is an $\vare$-coboundary expander then it is an $\vare^1$-topological expander.
\end{thm}

Now, combining Theorem 3.5 with Theorem 3.4, one deduces that $\D^{(d)}_n$ are topological expanders as mentioned in \S 3.1.

But these  are of unbounded degree.  Naturally,  as the finite quotients of the high rank Bruhat-Tits building are spectral  and geometric expanders, one tends to believe that they are also topological expanders.  This is still an open problem.  Let us say right away that in general these quotients (and even the Ramanujan complexes) are not coboundary expanders.  As was explained in \cite{KKL16} for many of the lattices $\Ga$ in simple $p$-adic Lie groups, $H^1(\Ga\setminus\mathcal{B}, \bbf_2) \neq \{ 0\}$ since it is equal to $H^1(\Ga, \bbf_2) = \Ga/[\Ga, \Ga]\Ga^2 $ (since $\mathcal{B} $ is contractible) and the latter is often non-zero.  Thus $\mathscr{h}_i (\Ga\setminus\mathcal{B}) = 0$ and $\Ga\setminus\mathcal{B} $ is not a coboundary expanders.

Still, one can overcome this difficulty.  For this we need another definition:
\begin{definition} A $d$-dimensional complex $X$ is called $\vare$\emph{-cosystolic expander}, if for every $i = 0,\dots, d-1$, one has $\nu_i (X) \ge \vare$ and $\mu_i (X) \ge \vare $ when:
\[ \nu_i (X) = \min\limits_{f \in C^i\setminus Z^i} \; \frac{\| \d_i(f)\|}{\|\lceil  f\rceil\|}\]
where $
\lceil f \rceil = f+ Z^i$ and
\[ \| \lceil f \rceil\| = \min\{ \| g \|\, \Big| \, g \in \lceil f \rceil \} \] and
\[ \mu_i = \min\limits_{f \in Z^i\setminus B^i} \| f\|.\]
\end{definition}

For later use, let us denote $\mu(X) = \min\mu_i(X)$
and $\nu(X) = \min \nu_i(X)$.  So, $X$ is $\vare$-cosystolic expansion if $\mu(X) \ge \vare $ and $\nu (X) \ge \vare$.
 So, $X$ is ``$\vare$-\emph{cocycle expander}";  it may not be coboundary expander if $H^i\neq \{ 0 \}$ (for some $i = 0,\dots, d-1$) but at least every representative of a non-trivial cohomology class is ``large".

An extension of Gromov's Theorem 3.5 is given in \cite{DKW16}:

\begin{thm}
Cosystolic expanders are topological expanders.
\end{thm}

It is natural to conjecture that the Ramanujan complexes and more generally the quotients of the high rank Bruhat-Tits buildings, while not
coboundary expanders,  in general,  are still cosystolic expanders.  But also this is open.  What is known is a somewhat weaker result which still suffices to answer,  in the affirmative, Gromov's question on the existence of bounded degree topological expanders.  The following theorem was proved by Kaufman-Kazhdan-Lubotzky \cite{KKL16} for $d \le 3$ and by Evra and Kaufman \cite{EvKa16}  for general $d$.

\begin{thm}
Fix $2 \le d \in \bbn$, then there exists $\vare = \vare (d) > 0$ and $q_0 = q_0(d)$ such that if $K$ is a local non-Archimedean field of fixed  residue degree  $q > q_0$ and $G = \utilde G(K)$ with $\utilde G$ simple $K$-group of $K$-rank $d$, then the $(d-1)$-skeletons $Y$ of the finite ($d$-dimensional) quotients $X$ of the Bruhat-Tits building $\mathcal{B} = \mathcal{B}(G)$ form a family of bounded degree $(d-1)$-dimensional $\vare$-cosystolic expanders.
\end{thm}

As this Theorem holds for every $d$, it solves Gromov's problem, but in a somewhat unexpected way.  We do believe that $X$ in the theorem are also cosystolic expanders and not just $Y$.

  Evra and Kaufman in \cite{EvKa16}, give a quite general combinatorial criterion to deduce a result like Theorem 3.8.  They prove that if $X$ is a $d$-dimensional complex of bounded degree all of whose proper links (i.e.\  $\ell k_X (F)$ for every face $F \neq \emptyset$) are coboundary expanders, and all the underlying graphs of all the links (including $\ell k_X(\emptyset) = X$) are ``very good" expander graphs, then the $(d-1)$-skeleton of $X$ is a cosystolic expander.  The reader is referred to \cite{EvKa16} for the exact quantitative formulation.  It is in spirit an ``$\bbf_2$-version" of Garland's local to global method described in \S 2.2.  It will be interesting to strengthen this result to the same level as Garland's, i.e., to assume only that the proper links are coboundary expanders and connected and $X$ is connected.  It will be even more interesting if one could deduce (even with the current hypothesis) that $X$ itself is  a cosystolic expander.  This will show that the $d$-dimensional Ramanujan complexes are topological expanders and not merely their $(d-1)$ skeletons as we now know.

The issue discussed in this section is only the tip of the iceberg. There are many more interesting problems (see \cite{Gro10}, \cite{GrGu12}) e.g.\ every $d$-dimensional complex can be embedded in $(2d+1)$-dimensional Euclidean space, but only some can be embedded in $2d$. Prove that high dimensional expanders (in some or any of the definitions) can not.

\section{Random simplicial complexes}

As mentioned briefly above, the easiest way to prove existence of bounded degree expander graphs is by random methods.  One may hope that this can be extended to the higher dimensional case of $d$-dimensional simplicial complexes.  But, here the problem is much more difficult.  In fact, as of now, there is no known ``random model" for $d$-dimensional simplicial complexes \emph{of bounded degree} (in the strong sense - see below) which gives high dimensional topological expanders.  This is surprising as the existence of such topological expanders is known by now by (\cite{KKL16}, \cite{EvKa16}) as  was explained in \S 3.  One may start to wonder if such a model exists at all, or maybe topological bounded degree expanders of high dimension are very rare objects.  Perhaps there is a kind of rigidity phenomenon analogue to what is well known by now in  Lie theory and locally symmetric spaces:  While there are many different  Riemann surfaces (parameterized by Teichm\"uller spaces), the higher dimensional case is completely different and rigidity results say that there are  ``very few" and mainly the ones coming from arithmetic lattices.

Let us now leave aside such a speculation and  give a brief background and a short account of the known results:

Erd\"os and R\'enyi initiated the study of random graphs in their seminal paper \cite{ER60}.  Their model is the following:  Let $n \in \bbn$ and $p \in [0, 1]$, the random model $X(n, p)$ is the graph $X$ with vertex set $[n] = \{ 1, \dots, n\}$ and where for every $1 \le i\neq  j \le n$, the edge $\{ i, j\}$ is in $X$ with probability $p$, independently of all other edges.  They  then study the properties of such graphs when $n \to \infty$ (and $p$ can be changed with $n$).  For example, their first famous result is that $p_0 = \frac{\log n}{n}$ is a threshold for the connectedness of $X \in X(n, p)$.  Namely,  for every $\vare > 0$,  if $p \le (1-\vare) \frac{\log n}{n}$ then almost surely (a.s.) such an $X$ is not connected, i.e.\  Prob. $(X \in X(n,p): X$ connected) $\underset{n\to\infty}{\longrightarrow} 0$.   On the other hand, if $p \ge (1 + \vare) \frac{\log n}{n}$ then $X$ is a.s. connected.

Why is $p = \frac{\log n}{n}$  the threshold?  Recall the ``coupon collector problem" which asserts that if elements of $[m] = \{ 1, \dots, m\}$ are chosen independently at random with repetition, it will take $t = m \log n$ steps  to choose them all.  In our process $p{n\choose 2}$ edges are chosen, and hence $2p{n\choose 2}$ vertices.  Now, if $p < \frac{\log n}{n}$ then less than $2\frac{\log n}{n} \frac{n^2}{2} = n \log n$ vertices are chosen.  So w.h.p there is an isolated vertex!  The amazing point in the Erd\"os-R\'enyi result is the fact that once we cross the threshold, not only are there no isolated vertices, but the graph is connected, and, in fact, even an expander.

This was the starting point of a very elaborate (and very important) theory of random graphs studying more and more delicate  properties of such $X \in X(n, p)$.

In \cite{LM06} Linial and Meshulam initiated  such a theory for 2-dimensional simplicial complexes.  A theory which shortly afterward was extended in \cite{MW09} to the general $d$-dimensional case.  The model studied $X^d(n,p)$ (nowadays called the Linial-Meshulam model for random  $d$-dimensional simplicial complexes) is the following:  $X \in X^d(n, p)$
 is a $d$-dimensional complex with $[n]$ as the set of vertices, $X$ contains the full $(d-1)$-skeleton, i.e., every subset of $[n]$ of size at most $d$ is in $X$ and a subset of size $d+1$ is in $X$ with probability $p$, independently of the other $d$-cells.  So $X^1 (n, p)$ is exactly the Erd\"os-R\'enyi model.  Now, for $d \ge 2$,  such an $X$  is always connected.  But, note that $X\in X^1(n, p)$ is connected if and only if $H^0 (X, \bbf_2) = \{ 0 \}$, so Linial, Meshulam and Wallach study for $d \ge 2$  and $X \in X^d(n, p)$, when $H^{d-1} (X, \bbf_2) = \{ 0 \}$ and proved the following far reaching generalization of the  Erd\"os-R\'enyi theorem.

 \begin{thm}[\cite{LM06} for $d=2$, \cite{MW09} for all $d$] The threshold for the homological connectivity, i.e.\  the vanishing of $H^{d-1} (X, \bbf_2)$ for $X \in X^d (n, p)$ is $p_0 = \frac{d\log n}{n}$.
 \end{thm}

The heuristic here for $\frac{d\log n}{n}$ is similar to the one above:  The process picks $p{n \choose d+1}$ $d$-cells and hence $(d+1) p {n\choose d+1} \; (d-1)$-cells.  So, if $p < \frac{d\log n}{n}$ less than $(d+1) \frac{d\log n}{n} {d \choose d+1} \approx {n\choose d} \log ({n \choose d}) \; (d-1)$-cells are chosen and so there is a $(d-1)$-cell $\tau$ with no $d$-cell containing it.  Hence the coboundary of $e_\tau$ - the characteristic function of $\tau$ - is zero, i.e. $e_\tau \in Z^{d-1} (X, \bbf_2)$.  On the other hand $e_\tau$ is not a coboundary  (note that in the complete $d$-dimensional complex $\delta(e_\tau) \neq 0$, so it is not even a cocycle)  and hence $H^{d-1} (X, \bbf_2) \neq \{ 0 \}$.  Again the interesting aspect of the Linial-Meshulam-Wallach result is that once the threshold is passed, not only does $H^{d-1} (X, \bbf_2)$ vanish, but $X$ is even a coboundary expander.

 A nice theory of random complexes has started to emerge (see \cite{Ka14} and the references therein).  As our main interest here is in expanders, we refer mainly to \cite{LM06}, \cite{MW09}, and \cite{DK12}, noting that  the results there  imply  (just like in the case of graphs) that  above the threshold the complexes are not only  homologically connected but also coboundary expanders and therefore topological expanders.
 The papers \cite{GuWa16} and \cite{KR16}  bring in a very detailed study of the spectrum of the higher dimensional Laplacians $\D_i$ action on $C^i (X,\bbr)$ for random $X$.

 But our main interest is in bounded degree complexes.  Recall that Bollobas \cite{Bo82} and others (see \cite{Wo99} for a comprehensive survey) have developed a  theory of random $k$-regular  graphs (for a fixed $k$) which also got a lot of attention.  This model, for $k\ge 3$, gives almost surely expander graphs of bounded degree.

 One would like to have such a model for $d$-dimensional complexes.  But first, what do we mean by bounded degree?  There are two natural meanings in the literature, which coincide  for $d = 1$.

 \begin{defn} A pure $d$-dimensional simplicial complex $X$ is of degree at most $k$ if every vertex of it is contained in at most $k$ cells of dimension $d$ (and so in at most $2^d\cdot k$ cells of any dimension). It is of \emph{upper-degree} at most $k$, if every face of dimension $d-1$ is contained in at most $k$ cells of  dimension $d$.
 \end{defn}

 A natural model of bounded degree simplicial complexes $Y^d(n, k)$ is given in \cite{FGLNP12}:  Assume, for simplicity, that $(d+1) | n$ and take a random partition of $[n]$ into $\frac{n}{d+1}$ subsets each of size $d+1$.  Choose independently $k$ such partitions and let $Y$ be the simplicial complex obtained by taking its cells to be all these $k$ $\frac{n}{d+1}$ subsets as well as all their subsets.  The case $d = 1$ boils down to the standard model of Bollobas.

 \begin{thm}[\cite{FGLNP12}]  For every fixed $d \in \bbn, \, \exists k_0 = k_0(d)$, such that for every $k\ge k_0$, a complex $Y \in Y^d(n, k)$ is almost surely $d$-dimensional geometric expander.
\end{thm}

This theorem is very promising at first sight, but unfortunately, $Y \in Y^d (n, k)$ is typically neither coboundary expander nor topological expander.  To visualize this think about the $d=2$ case:  When $k$ is fixed and $n$ very large, for a typical $Y \in Y^2(n, k)$, every edge of $Y$ is contained in at most one  triangle.  So, homotopically $Y$ looks more like a graph and one can map it into $\bbr^2$ with only small size overlapping points.

So, altogether, this is a nice model which certainly deserves further study (e.g.\ what is the threshold for $k_0 = k_0(d)$ in Theorem 4.3?) but it will not give us the stronger versions of expansion (topological, cosystolic, coboundary etc.).  As hinted at the beginning of this section,  it is still a major open problem to find a random model (if such at all exists) of $d$-dimensional bounded degree simplicial complexes which will give, say, topological expanders.

The situation with bounded \emph{upper} degree is better:  In \cite{LM15} Lubotzky and Meshulam gave  a model for 2-dimensional complexes of bounded upper degree
  (using the theory of Latin squares) and it was shown to produce coboundary expanders (and so also  topological expanders).  This was  generalized to all $d$ by Lubotzky-Luria-Rosenthal \cite{LLR15}, with a slight twist of the construction, replacing the Latin squares by Steiner systems and using the recent breakthrough of Keevash \cite{Kee14} on existence of designs.    Let us briefly describe the general model $W^d(n, k)$.

 Let $r \le q \le n$ be  natural numbers and $\lambda \in \bbn$.  An $(n, q, r, \lambda)$-\emph{design} is a collection $S$ of $q$-element subsets of $[n]$ such that each $r$-element  subset of $[n]$ is contained in exactly $\lambda$ elements of $S$.   Given $n, d \in\bbn$, an $(n, d)$-Steiner system is an $(n, d+1, d, 1)$-design, namely, a collection~$S$  of subsets  of size $d+1$ of $[n]$, such that each set of size $d$ is contained in exactly one element of $S$.  Using the terminology of simplicial complexes, an $(n, d)$-Steiner system can be considered as a $d$-dimensional simplicial complex of upper degree one.  Recently, in a groundbreaking paper \cite{Kee14}, Peter Keevash gave a randomized construction of Steiner systems for any fixed $d$ and large enough $n$ satisfying certain necessary divisibility conditions (which hold for infinitely many $n \in\bbn$).  From now on, we will assume that given a fixed $d \in \bbn$, the value of $n$ satisfies the divisibility condition from Keevash's theorem.

 Keevash's construction of Steiner systems is based on randomized algorithm which has two stages.  We will explicitly describe the first stage and use the second stage as a black box.

 Given a set of $d$-cells $A \subseteq {[n]\choose d+1}$, we call a $d$-cell $\tau$ \emph{legal with respect to } $A$ if there is no common $(d-1)$-cell in $\tau$ and in any cell in $A$.  Non-legal cells are also called forbidden cells.

 In the first stage of Keevash's construction, also known as the greedy stage, one selects a sequence of $d$-cells according to the following procedure.  In the first step, a $d$-cell is chosen uniformly at random from ${[n]\choose d+1}$.  Next, at each step a legal $d$-cell (with respect to the set of $d$-cells chosen so far) is chosen  uniformly at random and is added to the collection of previously chosen $d$-cells.  If no such $d$-cell exists the algorithm aborts.  The procedure stops when the number of $(d-1)$-cells which do not belong to the boundary of the chosen $d$-cells is at most $n^{d-\delta_0}$ for some fixed $\delta_0 > 0$ which only depends on $d$.  In particular, if the algorithm does not abort the number of steps is at least $({n \choose d} - n^{d-\delta_0}) / (d+1) \ge n^d / (2(d+1)!)$.

In the second stage, Keevash gives a randomized algorithm that adds additional $d$-cells in order to cover the remaining $(d-1)$-cells that are not contained in any of the $d$-cells chosen in the greedy stage.  We do not need to go into the details of this algorithm.  The important thing for us is that with high probability the algorithm produces an $(n, d)$-Steiner system.

Fix $k \in \bbn$ and let $S_1, \dots, S_k$ be $k$ independent copies of $(n, d)$-Steiner systems chosen according to the above construction, and let $X$ be the $d$-dimensional  simplicial complex whose $d$-cells are $\mathop{\bigcup}\limits^{k}_{i = 1} S_i$, so $X$ contains the complete $(d-1)$-skeleton and it is of upper degree at most $k$.

We can now state the main result of \cite{LLR15}:
\begin{thm} Fix $d \in \bbn$, there exists $k_0= k_0(d)$ and $\vare = \vare(d)$, such that for every $k\ge k_0$, a random complex $W \in W^d(n, k)$ is almost surely an $\vare$-coboundary expander, and hence also a topological expander.
\end{thm}

It will be of great interest to study various other properties of this model.  For example, find the threshold for $k_0(d)$ (the estimates obtained from \cite{LLR15} are huge and it will be very interesting to give more realistic upper bound, note that for $d = 1, k_0(d) = 3)$.  Another interesting problem is to study $\pi_1 (W)$-the fundamental group of $W$;  when is it hyperbolic? has property $(T)$? trivial?  The model $W$ behaves w.r.t. the model $X$ as Bollobas' model w.r.t. Erd\"os-R\'enyi, and this suggests many further directions of research on these bounded upper degree complexes.

\section{High dimensional expanders and computer science}

In recent years high dimensional expanders have captured the interest of computer scientists and various connections and applications have popped up. Most of these works are in their infancy. We will give here only a few short pointers on these developments, with the hope and expectation that the future will bring much more.

{\bf Probabilistically Checkable Proofs:} The PCP theorem, proven in the early 90's (cf. \cite{AS,ALMSS}), is a cornerstone of modern computational complexity theory stating that proofs can be written in a robust locally-testable format. PCPs are related to many areas within theoretical computer science  ranging from hardness of approximation to delegation and efficient cloud computing.

The basic PCP theorem can be proven using an expander-graph-based construction \cite{Din07}. For stronger PCPs, e.g. with unique constraints, or shorter proof length, or with lower soundness error, stronger forms of expansion seem to be needed, in particular high dimensional expansion might play a pivotal role. Dinur and Kaufman \cite{DK17} explore replacing the standard direct product construction (also known as parallel repetition \cite{Raz98}) by a much more efficient bounded-degree high dimensional expanders as constructed in \cite{LSV05}. Direct products are ubiquitous in complexity, especially as a useful hardness amplification construction, and bounded-degree high dimensional expanders may potentially be useful in many of those settings.

{\bf Locally testable codes:} LTCs are an information-theoretic analog of PCPs. These error correcting codes have the additional property that it is possible to locally test whether or not a received word is close to being a codeword. Unlike many problems in coding theory, this is a property that random codes {\em do not} have. This makes it even more challenging to settle the problem whether LTCs can have both linear rate and distance. The current best construction comes from a PCP and its rate is inverse poly-logarithmic \cite{B-SS08,Din07}.
High dimensional expanders naturally yield locally testable codes, whose parameters are unfortunately sub-optimal.

{\bf Property testing:} The central paradigm in property testing is the interplay between local views of an object and its global properties. The object can be a codeword, an NP-proof, or simply a graph. This theory generalizes both PCPs and LTCs and has significant practical applications. It was an unexpected discovery that high dimensional expanders (and especially the cohomological/coboundary expanders mentioned above) fit very naturally into this theory \cite{KL14}. Specifically, theorems about high dimensional expanders readily translate to results on property testing.

{\bf Quantum computation and quantum error correcting codes:} Sipser and Spielman \cite{SS96} showed how extremely good expander graphs yield excellent LDPC error-correcting codes. However, the existence of LDPC quantum error-correcting codes (even inexplicitly) remains a major open problem. Recent work by Guth and Lubotzky \cite{GuLu14} is a  step in this direction, which is related  to our topic: Every simplicial complex gives a ``homological error correcting code" (see \cite{BM07}, \cite{Ze09}) but in general they are of poor quality.  High dimensional coboundary expanders are related to local testability of codes (see \cite{AE15}).

Another  basic problem in quantum computation seeks a finite {\em universal set of quantum gates} that can efficiently generate an arbitrary unitary matrix in $U(n)$ to desired accuracy. This is solved by Kitaev and Solovay's classical algorithm, but non-optimally. The generators of Lubotzky-Phillips-Sarnak's Ramanujan graphs \cite{LPS86} fare better, but come with no efficient generative algorithm. Following the breakthrough of Ross-Selinger \cite{RS16}, the case $n=2$ is essentially solved in a recent work by Sarnak and Parzanchevski \cite{PS17} who came up with optimal (a.k.a.\ {\em golden}) gates and an explicit generative algorithm based on Ramanujan graphs. In ongoing work they use higher dimensional Ramanujan complexes to find such ``golden gates" for higher $n$ \cite{PS18}.

\medskip
\noindent
Acknowledgments: The author is indebted to S. Evra, G. Kalai, N. Linial,   O. Parzanchevski and R. Rosenthal  for remarks on an earlier version of this paper which greatly improved it. In addition, thanks are due to  ERC, NSF and BSF their support.



\bigskip
\centerline{Instititue of Mathematics}
\centerline{Hebrew University}
\centerline {Jerusalem 9190401}
\centerline{Israel}
\centerline{alex.lubotzky@math.huji.ac.il}
\end{document}